\documentclass[11pt]{article}
\usepackage{graphicx}
\usepackage{amsfonts}
\usepackage{bbm}
\usepackage{mathrsfs}
 \usepackage{mathrsfs,amsmath,amssymb}
\usepackage[dvips]{color}
 \usepackage{amssymb}
 \usepackage{amsbsy}
 \usepackage[dvips]{color}
\allowdisplaybreaks
 \newtheorem{condition}{\noindent\mbox{Condition}}
 \newtheorem{remark}{\noindent\mbox{Remark}}
 
 \newtheorem{lemma}{\noindent\mbox{Lemma}}
 \newtheorem{theorem}{\noindent\mbox{Theorem}}

 \def\bq{\begin{equation}}
 \def\eq{\end{equation}}
 \def\eqn{\end{eqnarray}}
 \def\bqn{\begin{eqnarray}}
 \def\proof{\noindent{\it Proof.~~}}
 \def\qed{\hfill$\Box$\medskip}
 \def\rto{\rightarrow\infty}
 \def\z{\left}
 \def\y{\right}

 \def\mf{\mathbf}
 \def\bs{\boldsymbol}

\begin{document}
\noindent{\Large\bf  On total progeny of multitype Galton-Watson process and the first passage time of random walk with bounded jumps} 

\noindent{
 Huaming Wang
}

\vspace{0.2 true cm}
\noindent{\footnotesize\rm
 Department of Mathematics, Anhui Normal University, Wuhu 241000, P. R. China;

 \noindent Email:hmking@mail.ahnu.edu.cn

\vspace{4mm}}%

\begin{center}
\begin{minipage}[c]{12cm}
\begin{center}\textbf{Abstract}\end{center}
In this paper, we first form a method to calculate the probability generating function of the total progeny of multitype branching process. As examples, we calculate probability generating function of the total progeny of the multitype branching processes within random walk which could stay at its position and (2-1) random walk. Consequently, we could give the probability generating functions and the distributions of the hitting time of corresponding  random walks.
\vspace{0.2cm}

\mbox{}\textbf{Keywords:}\quad Multitype branching process, total progeny, random walk.\\
\mbox{}\textbf{MSC 2010}:  Primary 60J80;
secondary 60G50.
\end{minipage}
\end{center}


\section{Introduction}
\subsection{Motivation}
Let $\{X_n\}$ be a simple random walk, that is, $\{X_n\}$ is a Markov chain with initial value $0$ and transitional probability $P(X_{n+1}=X_n+1|X_n,...,X_0)=p=1-P(X_{n+1}=X_n+1|X_n,...,X_0),$ where $0<p<1.$ Define $T=\inf[n>0:X_n=1]$ the hitting time of position $1.$ One follows, for example from reflection principle, that \begin{equation}\label{disnt}
  P(T=2n+1)=\frac{1}{2n+1}\left(\begin{array}{c}
                                           2n+1 \\
                                           n+1
                                         \end{array}\right)p^{n+1}(1-p)^n, \ n\ge 0.
\end{equation}
 Also one has that the probability generating function of $T,$ \begin{equation}\label{pgfnt}
   \zeta(s):=E(s^T)=\frac{1-\sqrt{1-4pqs^2}}{2qs},\ |s|<1.
 \end{equation}
For the specific
calculation of (\ref{disnt}) and (\ref{pgfnt}), see for example,  Strook \cite{s05}.
Next, we consider some more general random walk.

\vspace{0.2cm}
\noindent I). {\it Random walk with stay}

Suppose that $\{X_n\}$ is a  random walk which could stay at its position with positive probability $0<r<1,$ that is, $P(X_{n+1}=X_n+1|X_n,...,X_0)=p, $ $P(X_{n+1}=X_n-1|X_n,...,X_0)=q$ and $P(X_{n+1}=X_n|X_n,...,X_0)=r$ where $p,q,r>0$ and  $p+q+r=1.$
We call  $\{X_n\}$ a {\it random walk with stay.}

\vspace{0.2cm}
\noindent II). {\it (2-1) random walk}

Suppose $q_1,q_2>0$ and $p>0$ are three numbers such that $q_1+q_2+p=1.$
Let $\{X_n\}$ be a Markov chain, starting from $0,$ with transition probabilities
\begin{equation*}\begin{split}
  &P(X_{n+1}=X_n-l|X_n,...,X_0)=q_l ,\ l=1,2;\\
  &P(X_{n+1}=X_n+1|X_n,...,X_0)=p.
\end{split}
\end{equation*}
We call such $\{X_n\}$ an {\it (2-1) random walk.}

\vspace{0.2cm}
\noindent III). {\it (L-R) random walk}

 Let $\Lambda:=\{-L,...,R\}/\{0\}$ where $L$ and $R$ are two positive integers and  let $p_l,\ l\in \Lambda$ be positive numbers such that $\sum_{l\in\Lambda}p_l=1.$
Suppose  $\{X_n\}$ is a Markov Chain  with initial value $0$ and transition probabilities $P(X_{n+1}=X_n+l|X_n,...,X_0)=p_l,\ l\in \Lambda.$
 We call such $\{X_n\}$ a random walk with bounded jumps or simply an {\it  (L-R) random walk.}

 \vspace{0.2cm}
 For the random walks described above, we also define
 $$T=\inf[n> 0: X_n\ge 1] \text{ and } \zeta(s)=E(s^T),\ |s|<1.$$
 Stimulated by (\ref{disnt}) and (\ref{pgfnt}),  one  natural question is,
 could we give the distribution of $T$ and calculate its probability generating
 function $\zeta(s)?$
  In general, for the random walks described above, it is hard to give the distribution and
  the probability generating function of $T.$ Indeed, in the literatures we are aware of,
  no such result was given.

  We note that in  \cite{w12, hw10,hw10b, hzh10}, the authors revealed the intrinsic branching
  structures within random walk with stay, (L-1) random walk, (1-R) random walk,
  and (L-R) random walk.  The branching structure for  simple random walk  was given in Kesten Kozlov and Spitzer \cite{kks75}. In the above mentioned literatures, it has been found that the first passage time $T$ could be expressed in terms of the total progeny of the branching process. For simple random walk, the branching process is of single type, while for the other cases, the corresponding branching process is of multitype.

For general random walk, though it's hard to tell the distribution of $T$
 there are good news. Since the first passage time could be expressed in terms of the total progeny of branching process, once we give the distribution of the total progeny of multitype branching process, the distribution of $T$ follows.

For the single type Galton-Watson process, the probability generating function of the total progeny could be found in Dwass \cite{d68,d69} and also in Feller \cite{fe}.

In this paper, we first study the distribution of the total progeny of multitype branching process. We give a method to calculate the probability generating function of the total progeny. As two examples, we calculate probability generation functions of the total progeny of the multitype branching process within (2-1) random walk and random walk with stay. Consequently, we could calculate the probability generating function of $T$ explicitly. Especially, for the multitype branching process within the random walk with stay , in critical case, we show that, the tail probability of the total progeny decays as $\frac{C_1}{\sqrt{n}}$ for some constant $0<C_1<\infty$ when $n\rto.$ Thus, we show that the tail probability of the hitting time $T$ decays as $\frac{C_2}{\sqrt{n}}$  for some constant $0<C_2<\infty$ when $n\rto.$

We remark that for general (L-R) random walk, since the branching structure was revealed in Hong-Wang \cite{hw10b}, the probability generating function of $T$ follows similarly as that of (2-1) random walk and random walk with stay. But the branching structure involves a $(1+...+L)(1+...+R)$-type branching process. We do not contain such tedious calculation in this paper. 

\subsection{The  main results}

Let $L$ be a positive integer. Suppose that $p_i(\cdot),\ i=1,...,L$ are probability measures on $\mathbb Z_+^L,$ with $\mathbb Z_+=\{0,1,2,...\}.$ Let $\{Z_n\}_{n=0}^\infty$ be an L-type Galton-Watson process with the offspring distributions $P(Z_{n+1}=\z(n^{(1)},...,n^{(L)}\y)|Z_{n}=\mathbf e_i)=p_{i}\z(n^{(1)},...,n^{(L)}\y),\ i=1,2,...,L,$ where $\mathbf e_i\in \mathbb Z_+^L$ with the $i$-th component $1$ and all the others $0.$
For $1\le i\le L$ define
\begin{equation*}\begin{split}
  \phi^{(i)}&(s^{(1)},...,s^{(L)})=E\z(\z(s^{(1)}\y)^{Z_{1}^{(1)}}\cdots \z(s^{(L)}\y)^{Z_{1}^{(L)}}|Z_0=\mathbf e_i\y)\\
  &=\sum_{n^{(1)},...,n^{(L)}\ge0}p_i\z(n^{(1)},...,n^{(L)}\y)\z(s^{(1)}\y)^{n^{(1)}}\cdots \z(s^{(L)}\y)^{n^{(L)}},\ | s^{(l)}|<1,  1\le  l \le L,
\end{split}
\end{equation*}
 being the probability generating function of $Z_1$ given $Z_0=\mathbf e_i.$
Introduce $$Y_n=\sum_{i=0}^nZ_{i}$$ being the total progeny of the first $n$ generations  of $\{Z_n\}.$
Let $Y:=\lim_{n\rto}Y_n,$  the total progeny of the branching process, possibly  being $\infty.$
Define  generating function of $Y_n$
$$G_n^{(i)}(s^{(1)},...,s^{(L)})=E\z(\z(s^{(1)}\y)^{Y_n^{(1)}}\cdots \z(s^{(L)}\y)^{Y_n^{(L)}}|Z_0=\mathbf e_i\y), \ | s^{(l)}|<1,  1\le  l \le L.$$
Let $m_{i,j}=\frac{\partial\phi^{(i)}}{\partial s^{(j)}}.$
 Then it is known that $M:=(m_{i,j})\in \mathbb R^{L\times L}$ is the offspring matrix of $\{Z_n\}$.
Let  $\pi^{(i)}=P(Z_n=\mathbf 0 \text{ for some } n |Z_0=\mathbf e_i)$ being the extinction probabilities.
It will be convenient to introduce the vector notations
\begin{equation*}\begin{split}
  &\mathbf s=(s^{(1)},...,s^{(L)}),\\
  &\boldsymbol{\phi}(\mathbf s)=(\phi^{(1)}(\mathbf s),...,\phi^{(L)}(\mathbf s)),\\
  &\mathbf G_n(\mathbf s)=(G_n^{(1)}(\mathbf s),...,G_n^{(L)}(\mathbf s)),\\
  &\boldsymbol{\pi}=(\pi^{(1)},...,\pi^{(L)}),\\
  &\mathbf 1=(1,...,1).
\end{split}
\end{equation*}
For simplicity, we write $\mf G_1(\mf s)$ as $\mf G(\mf s).$
For two vectors $\mf s=(s^{(1)},...,s^{(L)})$ and $\mf n=(n^{(1)},...,n^{(L)}),$ we define $\mf s\mf n=(s^{(1)}n^{(1)},...,s^{(L)}n^{(L)}),$ $\mf s^{\mf n}=\z(s^{(1)}\y)^{n^{(1)}}\cdots \z(s^{(L)}\y)^{n^{(L)}}$ and
 the relation $\mf s\ll \mf n$  means that $ s^{(i)}<n^{(i)}, \text{ for } i=1,...,L.$
 The norm of a vector $\mf s$ is defined by $|\mf s|:=\max[s^{(i)}:1\le i\le L].$
Vectors are always in bold and  $v^{(i)}$ always means the $i$-th component of vector $\mf v.$ But as the branching process is involved, we do not use bold symbol. For example, $Z_n$ means the $n$-th generation and $Y$ is the total progeny of the branching process. They are all vectors.
\begin{condition}\label{con}
 Suppose that $\bs\phi(\mf s)$ are not linear functions of $s^{(1)},..., s^{(L)}$ and that $M^{n_0}\gg 0$ (all entry of $M^{n_0}$ are positive) for some integer $n_0\ge 1.$
\end{condition}
\begin{remark} Let $\sigma$ be the eigenvalue of $M$ with  largest magnitude.
 Under Condition 1, one follows from Frobenius theorem that $\sigma$ is real and positive.  Let $\bs\pi$ be the smallest solution of $\mf u=\bs\phi(\mf u).$ It's known that $\bs\pi$ is the extinction probability of the multitype branching process $\{Z_n\}$ and $\bs \pi=\mf 1$ or $\ll \mf 1$ according  as  $\sigma\le1,$ or $\sigma>1.$  \end{remark}

\begin{theorem}\label{gt}
  Suppose that Condition 1 holds. Then the limit $\bs\rho(s):=\lim_{n\rto}\mf G_n(\mf s)$  exists and for fixed $\mf s$, $\bs\rho(\mf s)$ is the unique solution of equation \begin{equation}\label{re}
   \mf u=\mf s \bs\phi(\mf u).
  \end{equation}
  Moreover, if $\sigma\le 1,$ then $\bs\rho(\mf s) $ is an honest probability generating function.
\end{theorem}
\begin{remark}
  Suppose $\sigma\le 1.$ Then the process $\{Z_n\}$ is extinct and $\bs\rho(\mf s)$ is a probability generating function.  As the limit of $\{\mf G_n\},$  $\bs\rho(\mf s)$ is the probability generating function of $Y,$ that is, $\rho^{(i)}(\mf s)=E({\mf s}^Y|Z_0=\mf e_i),\ |\mf s|<1.$
\end{remark}
As application of Theorem \ref{gt}, we study two special 2-type branching process which are connected with random walk with stay and (2-1) random walk respectively. 
\begin{theorem}\label{srw}
Let $\{Z_n\}_{n\ge 0}$ be a 2-type branching process with offspring distributions
\begin{eqnarray}
   &&P\left(Z_{n+1}=(a,b)|Z_n=\mf e_1\right)=\frac{(a+b)!}{a!b!}q^ar^bp,\ a,b\ge 0 \label{bra}\\
    && P\left(Z_{n+1}=\mf 0|Z_n=\mf e_2\right)=1,\label{brb}
   \end{eqnarray}
   where $p,q,r>0$ and $p+q+r=1.$
 Let $Y_n=\sum_{i=0}^nZ_i$ and $Y=\lim_{n\rto}Y_n.$ Then we have
Suppose that $q\le p.$ Then $P(Y<\infty)=1$ and the probability generating function of $Y$
\begin{equation}\label{pgfrw}
\bs\rho(\mf s)=(\rho^{(1)}(\mf s), \rho^{(2)}(\mf s))=\Big(\frac{1-rs^{(2)}-\sqrt{\z(1-rs^{(2)}\y)^2-4pqs^{(1)}}}{2q}, s^{(2)}\Big).
\end{equation}
Moreover, if $p=q=\frac{1-r}{2}$ and $P(Z_0=\mf e_1)=1$ then \begin{equation}\label{tl}
\lim_{n\rto}\sqrt{n}P(|Y|>n)=\frac{1}{\sqrt\pi}\sqrt\frac{{1+r}}{1-r}.
\end{equation}
At last, if we replace (\ref{brb}) by 
\begin{equation}\label{brc}P\left(Z_{n+1}=(a+1,b)|Z_n=\mf e_1\right)=\frac{(a+b)!}{a!b!}q^ar^bp, \ a,b \ge0,\end{equation}
then $\rho^{(1)}(\mf u)$ is the smallest real solution of 
$$r\frac{u^{(2)}}{u^{(1)}}\z(\rho^{(1)}\y)^3-q\z(\rho^{(1)}\y)^2-\rho^{(1)}+pu^{(1)}=0$$ and $\rho^{(2)}=\z(\rho^{(1)}\y)^2\frac{u^{(2)}}{u^{(1)}}.$
 \end{theorem}
\begin{remark}
  The branching process $\{Z_n\}$ determined by (\ref{bra}) and (\ref{brb}) corresponds to branching structure for random walk with stay and the one determined by (\ref{bra}) and (\ref{brc}) corresponds to branching structure for (2-1) random walk.
\end{remark}
Since the hitting time $T$ of random walk could be expressed in terms of the total progeny $Y$ of the branching process. We study the distribution and generation of random walk in the following theorems.

\begin{theorem}\label{st}
  Suppose $\{X_n\}$ is a random walk with stay  and that $q\le p.$ Let $T$ be its first passage time of position $1.$ Then
  $$E(u^T)=\frac{1-ru-\sqrt{\z(1-ru\y)^2-4pqu^2}}{2qu},\ |u|<1.$$ Moreover if $p=q=\frac{1-r}{2}$ then
  $$\lim_{n\rto}\sqrt nP(T\ge n)=\sqrt\frac{2}{\pi}\frac{1}{\sqrt{1-r}}.$$
\end{theorem}

\begin{theorem}\label{l1}
  Suppose that $\{X_n\}$ is a (2-1) random walk and   that $p-q_1-2q_2\ge0.$ Let $h(u)$ be the probability generating function of the first passage time $T.$ Then $g(s):=sh(s)$ is the smallest real solution of equation
  \begin{equation}\label{te}q_2s^{-1}g^3+q_1g^2-g+ps=0.\end{equation}
\end{theorem}

\begin{remark}
  Results for random walk with stay  look more satisfied. We give explicitly the probability generating function of $T,$ and calculate its tail probability. Indeed, in the proof of the theorem, we do give its exact distribution $P(T=n).$
  As far as (2-1) random walk is concerned, we could also solve equation (\ref{te}) to get the probability generating function $h(s)$ of $T$ and expand $h(s)$ to give its exact distribution. But since the root is complicated, we omit such calculation.

  We mention that for general (L-R) random walk, in a similar method one could use the branching structure to give the probability generating function of $T.$ Of course, in this general case, the root for equation $\mf u=\bs\phi(\mf u)$ will be complicated.
\end{remark}
The remainder part of the paper is arranged as follows. Section \ref{bstr} is a short introduction of the branching structure within random walk with stay  and (2-1) random walk. Section 3 contains the proofs of the above theorems.
\section{Branching structures within random walks}\label{bstr}
In this section, we introduce the branching structure within random walk with stay and (2-1) random walk. Firstly,
consider random walk with stay.
Recall that  $T=\inf[n>0:X_n=1]$ is the hitting time of position $1.$
Let
\begin{equation}\label{ua}
   \z(U_1^{(1)},U_1^{(2)}\y)=(1,0) \text{ and for }i\le0,
\end{equation}
and for $i\le 0$ define
\begin{equation}\label{ub}
  \begin{split}
    &U_i^{(1)}=\#\{0\le n<T_1: X_n=i,X_{n+1}=i-1\},\\
    & U_i^{(2)}=\#\{0\le n<T_1: X_n=i,X_{n+1}=i\}
  \end{split}
\end{equation}
counting the number of steps by the walk from $i$ to $i-1$ and the steps from $i$ to $i$ itself respectively.
Then it has been shown that $\z\{\z(U_n^{(1)},U_n^{(2)}\y)\y\}_{n\le 1}$ is a 2-type branching process. We summarize the branching structure of random walk with stay $\{X_n\}$ in the following theorem. For details of the proof, we refer the reader to Wang \cite{w12}  and Zeitouni \cite{ze04}.

\noindent {\textbf{Theorem A:}} {\it Let $\{X_n\}$ be a random walk with stay. If $q\le p,$ then  $\z\{\z(U_n^{(1)},U_{n}^{(2)}\y)\y\}_{n\le 1}$  defined in (\ref{ua}) and (\ref{ub}) forms a  $2$-type branching process. Its offspring distributions are
 \begin{equation*}\label{micha}\begin{split}
   &P\left(\z(U_i^{(1)},U_i^{(2)}\y)=(a,b)|\z(U_{i+1}^{(1)},U_{i+1}^{(2)}\y)=(1,0)\right)=\frac{(a+b)!}{a!b!}q^ar^bp,\ a,b\ge0, \\
    & P\left(\z(U_i^{(1)},U_i^{(2)}\y)=(0,0)|\z(U_{i+1}^{(1)},U_{i+1}^{(2)}\y)=(0,1)\right)=1.
   \end{split}
   \end{equation*}
 Moreover  the hitting time $T$ could be expressed by the branching process  as  	
	\begin{equation*}\label{ut}
T=1+\sum_{i\le0}2U_i^{(1)}+U_{i}^{(2)}.
\end{equation*}	}
Next, consider (2-1) random walk.  Let $T=\inf[n>0:X_n=1].$
Define, for $-\infty<i\le0,$
 \begin{equation}\label{uc}
   \begin{split}
     &U^{(1)}_{i}=\#\{0<k<T_1:X_{k-1}>i,X_k=i\}\\
 &U^{(2)}_{i}=\#\{0<k<T_1:X_{k-1}>i,X_k=i-1\}
   \end{split}
 \end{equation}
and set \begin{equation}\label{ud}
  \z(U_1^{(1)},U_1^{(2)}\y)=(1,0).
\end{equation}
Then we have the following theorem, whose proof could be find in  Hong and Wang \cite{hw10}.

\noindent{\bf Theorem B:} {\it
Let $\{X_n\}$ is a (2-1) random walk. Suppose that $E(X_1)=p-q_1-2q_2\ge0.$ Then
$\z\{\z(U_i^{(1)},U_i^{(2)}\y)\y\}_{i\le 1}$ defined in (\ref{uc}) and (\ref{ud}) forms a 2-type branching process with offspring distributions
\begin{equation*}\label{michb}\begin{split}
   &P\left((U_{i-1}^{(1)},U_{i-1}^{(2)})=(a,b)\big|(U_{i}^{(1)},U_{i}^{(2)})=(1,0)\right)=\frac{(a+b)!}{a!b!}q_1^aq_2^bp, \\
    & P\left((U_{i-1}^{(1)},U_{i-1}^{(2)})=(a+1,b)\big|(U_{i}^{(1)},U_{i}^{(2)})=(0,1)\right)=\frac{(a+b)!}{a!b!}q_1^aq_2^bp, \ a,b\ge0 \end{split}
   \end{equation*}
and that $$T=1+\sum_{i\le0}2U_i^{(1)}+U_i^{(2)}.$$}

\section{Proofs}
\subsection{Proof of Theorem \ref{gt}}
We claim that for $n\ge 0,$ $\mf G_{n+1}(\mf s)=\mf s \bs\phi(\mf G_{n}(\mf s))$  with $\mf G_0(\mf s):=\mf s.$ Indeed,  note that $ G_1^{(i)}(\mf s)=E(\mf s^{Z_0+Z_1}|Z_0=\mf e_i)=s^{(i)}E(\mf s^{Z_1}|Z_0=\mf e_i)=s^{(i)}\phi^{(i)}(\mf s).$ Then one follows that $\mf G_1(\mf s)=\mf s\bs\phi(\mf s).$
Suppose that $\mf G_{m+1}(\mf s)=\mf s \bs \phi(\mf G_{m}(\mf s))$  for all $m< n.$ Let $f(\mf v|\mf u)=P(Z_1=\mf v|Z_0=\mf u).$ Then
\begin{eqnarray*} 
     &&G_{n+1}^{(i)}(\mf s)=E(\mf s^{Z_0+Z_1+...+Z_{n+1}}|Z_0=\mf e_i)\\
     &&\hspace{1.48cm}=s^{(i)}\sum_{\mf z_1,...\mf z_{n+1}}s^{\mf z_1+...+\mf z_n+\mf z_{n+1}}f(\mf z_1|\mf e_i)\cdots f(\mf z_{n}|\mf z_{n-1})f(\mf z_{n+1}|\mf z_n)\\
     &&\hspace{1.48cm}=s^{(i)}\sum_{\mf z_1,...\mf z_{n}}s^{\mf z_1+...+\mf z_n}f(\mf z_1|\mf e_i)\cdots f(\mf z_{n}|\mf z_{n-1})\sum_{z_{n+1}}\mf s^{\mf z_{n+1}}f(\mf z_{n+1}|\mf z_n)\\
     &&\hspace{1.48cm}=s^{(i)}\sum_{\mf z_1,...\mf z_{n}}s^{\mf z_1+...+\mf z_n}f(\mf z_1|\mf e_i)\cdots f(\mf z_{n}|\mf z_{n-1})(\bs\phi(\mf s))^{\mf z_n}\\
     &&\hspace{1.48cm}=s^{(i)}\sum_{\mf z_1,...\mf z_{n}}s^{\mf z_1+...+\mf z_{n-1}}f(\mf z_1|\mf e_i)\cdots f(\mf z_{n}|\mf z_{n-1})(\mf G_1(\mf s))^{\mf z_n}.
\end{eqnarray*}
Repeating the above procedure for $n$ times, it follows that $G_{n+1}^{(i)}(\mf s)=s^{(i)} \phi^{(i)}(\mf G_{n}(\mf s)).$  The claim is proved.

For $ \mf 0\ll \mf s \ll\mf 1,$ one easily sees that $\mf G_1(\mf s)=\mf s\bs \phi(\mf s)\ll \mf s=\mf G_0(\mf s).$ Then one follows from induction that $\{\mf G_n(s)\}$ is monotone decreasing in $n.$ Consequently, the limit $\bs\rho(\mf s):=\lim_{n\rto}\mf G_n(\mf s)$ exists and $\bs\rho(\mf s)$ satisfies equation 
\begin{equation}
\label{solu}\bs\rho(\mf s)=\mf s\bs\phi(\bs \rho(\mf s)).
\end{equation}
Next we show that for fixed $\mf 0\ll \mf s\ll \mf 1,$ equation (\ref{re}) has only one solution. Suppose $\bs \pi$ is the smallest nonnegative root of equation $$\mf u=\bs\phi(\mf u).$$ That is $\bs\pi$ is the extinction probability of the branching process. Then $\bs \pi\ll \mf1$ or $=\mf 1$ according as $\sigma\le 1$ or $>1.$ From equation (\ref{re}) we know that $u^{(1)}=s^{(1)}\phi^{(1)}\z(u^{(1)},u^{(2)},...,u^{(L)}\y).$ Noting that,
as a function of $u^{(1)},$ $y\z(u^{(1)}\y)=s^{(1)}\phi^{(1)}\z(u^{(1)},u^{(2)},...,u^{(L)}\y)$ is  strictly convex, thus $u^{(1)}=s^{(1)}\phi^{(1)}\z(u^{(1)},u^{(2)},...,u^{(L)}\y)$ has at most two roots. Since  $0<s^{(1)}\phi^{(1)}\z(0,u^{(2)},...,u^{(L)}\y)$ and $\pi^{(1)}>s^{(1)}\phi^{(1)}\z(\pi^{(1)},u^{(2)},...,u^{(L)}\y)$ one concludes that for fixed $u^{(2)},...,u^{(L)}$ there exists only one $0<u^{(1)}<\phi^{(1)}$ such that $u^{(1)}=s^{(1)}\phi^{(1)}\z(u^{(1)},u^{(2)},...,u^{(L)}\y).$ On the other hand, since $1>s^{(1)}\phi^{(1)}\z(1,u^{(2)},...,u^{(L)}\y),$ $u^{(1)}=s^{(1)}\phi^{(1)}\z(u^{(1)},u^{(2)},...,u^{(L)}\y)$ has no root in $[\pi^{(1)},1]$ Thus the solution of (\ref{re}) is unique.

Letting $\mf s=\mf 1$ in (\ref{solu}) one follows that $\bs \rho(\mf 1)=\bs\pi.$ Then $\bf\rho(\mf s)$ is an honest probability generating function if and only if $\sigma\le 1.$\qed

\subsection{ Proof of Theorem \ref{srw}} As beginning, we prove the first part. From the offspring distributions (\ref{bra}) and (\ref{brb}), one follows that $$\phi^{(1)}(s^{(1)},s^{(2)})=\sum_{a=0}^\infty\sum_{b=0}^\infty\frac{(a+b)!}{a!b!}q^ar^bp\z(s^{(1)}\y)^a\z(s^{(2)}\y)^b=\frac{p}{1-qs^{(1)}-rs^{(2)}}$$ and that $$\phi^{(2)}(s^{(1)},s^{(2)})\equiv 1.$$ Therefore

$$\bs\phi(\mf s)=(\frac{p}{1-qs^{(1)}-rs^{(2)}},1).$$

Next we form the probability generating function of the total progeny $Y.$ We note that Condition \ref{con} does not hold since the particles of the second type does not give offsprings. But this causes no problem to use Theorem \ref{gt} since the conditions on the mean offspring matrix $M$ are only used to ensure the extinction of the branching process. Indeed, one follows from (\ref{bra}) that 
$$P(Z_{n+1}^{(1)}=a|Z_{n}^{(1)}=1)=\z(\frac{q}{p+q}\y)^a\frac{q}{p+q}.$$ Thus,
if $q\le p,$ $\bs\pi=\mf 1,$ that is, the branching process $\{Z_n\}$ is extinct and $P(Y<\infty)=1.$

Let $\rho^{(i)}(\mf s)=E(\mf{s}^Y|Z_0=\mf{e}_i),\ i=1,2$ and $\bs\rho(\mf s)=(\rho^{(1)}(\bf s),\rho^{(2)}(\bf s)).$
Then One follows from Theorem \ref{gt} that $\bs\rho(\bf s)$ solves the equation
$$\bs\rho=\mf s \bs\phi(\bs\rho).$$ Consequently, $$\bs\rho(\mf s)=(\rho^{(1)}(\mf s), \rho^{(2)}(\mf s))=\Big(\frac{1-rs^{(2)}-\sqrt{\z(1-rs^{(2)}\y)^2-4pqs^{(1)}}}{2q}, s^{(2)}\Big).$$
The first part of Theorem \ref{srw} follows.

For the proof of the second part, letting $p=q=\frac{1-r}{2} $ in (\ref{pgfrw}), one follows that
if $P(Z_0=\mf e_1)=1,$ the probability generation function of $Y$ is
$$\rho^{(1)}(\mf s)=\frac{1-rs^{(2)}-\sqrt{\z(1-rs^{(2)}\y)^2-(1-r)^2s^{(1)}}}{1-r}.$$

Let $\beta(u):=E(u^{|Y|})=E(u^{Y_1}u^{Y_2}),\ |u|<1$ be the probability generating function of $|Y|.$ Then one sees that $$\beta(u)=\rho^{(1)}(u,u)=\frac{1-ru-\sqrt{(1-ru)^2-(1-r)^2u}}{1-r}.$$ Suppose that
$$\beta(u)=\beta_0+\beta_1u+\beta_2u^2+...\beta_nu^n+...$$ and define for $n\ge 0,$ $$\theta_n=\beta_{n+1}+\beta_{n+2}+....$$ Then $\theta_n=P(|Y|>n).$ Therefore, to prove (\ref{tl}), it is enough to show that \begin{equation}\label{tn}
 \lim_{n\rto}\sqrt n\theta_n=\frac{1}{\sqrt\pi}\sqrt\frac{{1+r}}{1-r}.
\end{equation}
For this purpose, define
$$\theta(u)=\theta_0+\theta_1u+\theta_2u^2+...\theta_nu^n+...,\ |u|<1.$$
One follows that \begin{equation}\label{the}\theta(u)=\frac{1-\beta(u)}{1-u}=\frac{1}{1-r}\sqrt\frac{{1-r^2u}}{{1-u}}-\frac{r}{1-r}.\end{equation}

In order to prove (\ref{tn}), we first prove the following lemma.
\begin{lemma}\label{vol}
  For $n\ge 0,$ let $a_n(x)=\frac{(2n-3)!!}{(2n)!!}x^n$ and $b_n=\frac{(2n-1)!!}{(2n)!!},$ where $|x|<1$ and $k!!=0$ for $k\le0.$  Then $$\lim_{n\rto}\sqrt n\sum_{k=0}^na_kb_{n-k}=\frac{1}{\sqrt\pi}(2-\sqrt{1-x}).$$
\end{lemma}
\proof
 Note that
\begin{equation}\label{tp}\begin{split}
   \sum_{k=0}^na_k&b_{n-k}\\
   &= \left.\begin{array}{l}
      \frac{(2n-1)!!}{(2n)!!}+\frac{(2n-3)!!}{(2n-2)!!}x+\frac{(2n-5)!!}{(2n-4)!!}\frac{x^2}{2} \\
      \hspace{1cm}+\frac{(2n-7)!!}{(2n-6)!!}\frac{3!!x^3}{6!!}+...+\frac{(2[\frac{n}{2}]-1)!!}{(2[\frac{n}{2}])!!}\frac{(2[\frac{n}{2}]-3)!!x^{[\frac{n}{2}]}}{(2[\frac{n}{2}])!!} \end{array}\right\}\cdots\cdots (\text{Part I})\\
      &\quad\quad+\sum_{k=0}^{[\frac{n}{2}]-1}\frac{(2n-2k-3)!!}{(2n-2k)!!}\frac{(2k-1)!!x^{n-k}}{(2k)!!}.\cdots\cdots\cdots\cdots \cdots\cdot \cdot(\text{Part II})
\end{split}
\end{equation}
Then using the fact $$\sqrt\frac{2}{\pi}\frac{1}{\sqrt{2n+1}}\le \frac{(2n-1)!!}{(2n)!!}<\sqrt\frac{2}{\pi}\frac{1}{\sqrt{2n}},$$ we have that (Part I) of (\ref{tp}) is smaller than
\begin{equation}\label{sqi} \sqrt\frac{2}{\pi}\Big(\frac{1}{\sqrt{2n}}+\frac{1}{\sqrt{2n-2}}\frac{x}{2}+\frac{1}{\sqrt{2n-4}}\frac{x^2}{4!!}+\frac{1}{\sqrt{2n-6}}\frac{3!!x^2}{6!!}+...+\frac{1}{\sqrt{2[\frac{n}{2}]}}\frac{(2[\frac{n}{2}-3])!!x^{[\frac{n}{2}]}}{(2[\frac{n}{2}])!!}\Big),
\end{equation}
and (Part II) of (\ref{tp}) is smaller than
\begin{equation*}\label{sqii}
\sqrt\frac{2}{\pi}\sum_{k=0}^{[\frac{n}{2}]-1}\frac{1}{\sqrt{(2n-2k)}}\frac{1}{2n-2k-1}\frac{(2k-1)!!x^{n-k}}{(2k)!!}. \end{equation*}
Denoting (\ref{sqi}) by $f(n),$ we have that
\begin{equation}\label{pli}
  \lim_{n\rto}\sqrt nf(n)=\frac{1}{\sqrt\pi}\sum_{n=0}^\infty\frac{(2n-3)!!x^n}{(2n)!!}=\frac{1}{\sqrt\pi}(2-\sqrt{1-x}).
\end{equation}
On the other hand, \begin{equation}\label{plii}
  \lim_{n\rto}\sqrt n\sqrt\frac{2}{\pi}\sum_{k=0}^{[\frac{n}{2}]-1}\frac{1}{\sqrt{(2n-2k)}}\frac{1}{2n-2k-1}\frac{(2k-1)!!x^{n-k}}{(2k)!!} =0.
\end{equation}
Then (\ref{tp}), (\ref{pli}) and (\ref{plii}) imply that $$\lim_{n\rto}\sqrt n\sum_{k=0}^na_kb_{n-k}\le \frac{1}{\sqrt\pi}(2-\sqrt{1-x}).$$
A similar argument yields that $$\lim_{n\rto}\sqrt n\sum_{k=0}^na_kb_{n-k}\ge \frac{1}{\sqrt\pi}(2-\sqrt{1-x}).$$
The lemma is proved.\qed

Next we continue the proof of Theorem \ref{srw}. Note that $$\frac{1}{\sqrt{1-u}}=1+\sum_{n=1}^\infty\frac{(2n-1)!!}{(2n)!!}u^n\text{ and } \sqrt{1-r^2u}=1-\sum_{n=1}^\infty\frac{(2n-3)!!}{(2n)!!}r^{2n}u^n.$$
Let $c_0=1, d_0=1$ and for $n\ge 1,$ let $c_n=-\frac{(2n-3)!!}{(2n)!!}r^{2n},$ $d_n=\frac{(2n-1)!!}{(2n)!!}.$
Then one has that
\begin{equation*}
 \sqrt \frac{1-r^2u}{1-u}=\sum_{n=0}^\infty \Big(\sum_{k=0}^nc_kd_{n-k}\Big)u^n.
\end{equation*}
Substituting to (\ref{the}), we have that $\theta_0=1$ and \begin{equation*}\begin{split}
\theta_n&=\frac{1}{1-r}\sum_{k=0}^{n}c_kd_{n-k}=\frac{1}{1-r}\Big(\frac{(2n-1)!!}{(2n)!!}-\sum_{k=1}^na_k(r^2)b_{n-k}\Big)\\
&=\frac{1}{1-r}\Big(\frac{2(2n-1)!!}{(2n)!!}-\sum_{k=0}^na_k(r^2)b_{n-k}\Big).
\end{split}\end{equation*}
Using Lemma \ref{vol}, one follows  that $$\lim_{n\rto}\sqrt n\theta_n=\frac{1}{1-r}\Big(\frac{2}{\sqrt\pi}-\frac{1}{\sqrt\pi}\left(2-\sqrt{1-r^2}\right)\Big)=\frac{1}{\sqrt\pi}\sqrt\frac{{1+r}}{1-r},$$
which proves (\ref{tn}).

The third part of Theorem \ref{srw} could be found in the proof of Theorem \ref{l1}. 
\qed
\subsection{Proof of Theorem \ref{st}}
Let $\{X_n\}$ be a random walk with stay. Comparing the branching process $\{Z_n\}$ in Theorem \ref{srw} and $\{(U_n^{(1)},U_n^{(2)})\}_{n\le1}$ in Theorem B (see Section \ref{bstr}), if $P(Z_0=\mf e_1)=1$ one follows that
 $T$ has the same distribution with $$1+\sum_{n=1}^\infty 2Z_{n}^{(1)}+Z_n^{(2)}=\sum_{n=0}^\infty 2Z_{n}^{(1)}+Z_n^{(2)}-1=2Y^{(1)}+Y^{(2)}-1.$$
Letting $\eta(u):=E(u^{T+1}),$ then $$\eta(u)=E(u^{2Y^{(1)}}u^{Y^{(2)}})=\rho^{(1)}(u^2,u)=\frac{1-ru-\sqrt{(1-ru)^2-(1-r)^2u^2}}{1-r}.$$
Define $\alpha_n=P(T\ge n)$ and let $\alpha(u)=\sum_{n=0}^\infty\alpha_nu^n.$ Similar as (\ref{the}) we have that $$\alpha(u)=\frac{1-\eta(u)}{1-u}=\frac{1}{1-r}\sqrt\frac{1-(2r-1)u}{1-u}-\frac{r}{1-r}.$$
Following the lines of the proof of the second part of Theorem \ref{srw}, replacing $r^2$ by $2r-1,$ we could prove that $$\lim_{n\rto}\sqrt n\alpha_n=\frac{1}{1-r}\Big(\frac{2}{\sqrt\pi}-\frac{1}{\sqrt\pi}\left(2-\sqrt{1-(2r-1)}\right)\Big)=\sqrt\frac{2}{\pi}\frac{1}{\sqrt{1-r}}.$$ \qed

\subsection {Proof of Theorem \ref{l1}}

We only give the idea of the proof since it is similar as that of Theorem \ref{st}.
Let $\{X_n\}$ be a (2-1) random walk. Suppose $\z\{\z(U_{i}^{(1)},U_{i}^{(2)}\y)\y\}_{i\le 1}$ is the process defined in Theorem B (see Section \ref{bstr}).
 Let $Z_n=(Z_n^{(1)},Z_n^{(2)})=(U_{-n+1}^{(1)},U_{-n+1}^{(2)})$ and $Y=\sum_{n=0}^\infty Z_n.$  Then  one follows from Section \ref{bstr}, Theorem B that $$T=2Y^{(1)}+Y^{(2)}-1.$$
 It is easy to calculate that  $$\bs\phi(\mf s)=(\phi^{(1)}(\mf s),\phi^{(2)}(\mf s))=\z(\frac{p}{1-q_1s^{(1)}-q_2s^{(2)}},\frac{ps^{(1)}}{1-q_1s^{(1)}-q_2s^{(2)}}\y).$$
Then one follows from Theorem \ref{gt} that $\bs\rho(\mf u)=(\rho^{(1)}(\mf u),\rho^{(2)}(\mf u))$ is the smallest real solution of $$\bs \rho=\mf u\bs\phi(\bs\rho).$$ That is, $(\rho^{(1)}(\mf u),\rho^{(2)}(\mf u))$ satisfies
$$(\rho^{(1)},\rho^{(2)})=\z( \frac{pu^{(1)}}{1-q_1\rho^{(1)}-q_2\rho^{(2)}},\frac{pu^{(2)}\rho^{(1)}}{1-q_1\rho^{(1)}-q_2\rho^{(2)}}\y).$$
Therefore, $\rho^{(1)}(\mf u)$ is the smallest real solution of $$q_2\frac{u^{(2)}}{u^{(1)}}\z(\rho^{(1)}\y)^3+q_1\z(\rho^{(1)}\y)^2-\rho^{(1)}+pu^{(1)}=0.$$
 Let \begin{equation*}\label{sol}\rho^{(1)}(u^{(1)},u^{(2)}):=F(q_2\frac{u^{(2)}}{u^{(1)}},q_1,-1,pu^{(1)})\end{equation*} be such a solution.

Noting that $sh(s)=E(s^{T+1})=E(s^{2Y^{(1)}}s^{Y^{(2)}})=\rho^{(1)}(s^2,s),$  one has that
$$g(s):=sh(s)=F(q_2s^{-1},q_1,-1,ps),$$ which is the smallest real solution of equation   $$q_2s^{-1}g^3+q_1g^2-g+ps=0.$$ Theorem \ref{l1} is proved.\qed

\noindent{\large{\bf Acknowledgements:}} The author wishes to thank Professor Wenming Hong, who advises him to consider this project. Thanks also extend to Miss Hongyan Sun for her discussion on the estimation of some specific series.



\begin{thebibliography}{99}
  \addtolength{\itemsep}{-0.5em}
        \bibitem{d68} Dwass, M. (1968) A theorem about infinitely divisible distribution. {\it Z. Wahrscheinlichkeitsth. 9, 287-289}
    \bibitem{d69} Dwass, M. (1969). The total progeny in a branching process and a related random walk. {\it J. Appl. Prob. 6, 682-628.}
         \bibitem{fe} Feller, W. (1968). An introduction to  probability theory and its application. {\it Vol. 2, 3rd edition. Wiley, New York.}
        \bibitem{hw10} Hong, W.M. and Wang, H.M. (2012).
Intrinsic branching structure within (L-1) random walk in random environment and its applications. {\it  To appear in Infinite Dimensional Analysis, Quantum Probability and Related Topics. }
\bibitem{hw10b} Hong, W. M. and Wang, H.M. (2010).  Intrinsic branching structure within random walk on $\mathbb{Z}.$ {\it arXiv:1012.0636.}
\bibitem{hzh10} Hong, W. M. and Zhang, L. (2010).
Branching structure for the transient (1,R)-random
walk in random environment and its applications. {\it  Infinite Dimensional Analysis, Quantum Probability and Related Topics, 13, pp. 589-618. }
    \bibitem{kks75} Kesten, H., Kozlov, M.V., and Spitzer, F. (1975). A limit law
for random walk in a random environment. {\it Compositio Math. 30,
pp. 145-168.}
\bibitem{s05}Strook, D.W. An introduction to Markov processes. (2005). {\it Springer verlag.  }
\bibitem{w12} Wang, H. M. (2012). Mean and variance of first passage time
of non-homogeneous random walk. {\it Front. Math. China, Vol. 7,  No.3,  551-559.}
\bibitem{ze04}  Zeitouni, O. (2004). Random walks in random environment. \textit{LNM 1837, J. Picard (Ed.), pp. 189-312, Springer-Verlag Berlin Heidelberg}.

\end{thebibliography}
  \end{document}